\documentclass[12pt]{amsart}

\newtheorem{thm}{Theorem}
\newtheorem{lem}{Lemma}

\newtheorem{cor}{Corollary}

\newtheorem*{thma}{Theorem A}
\newtheorem*{thmb}{Theorem B}
\newtheorem*{thm2'}{Theorem 2$^\prime$}

\newcommand{\vphi}{\varphi}

\newcommand{\Hess}{\operatorname{Hess}}
\newcommand{\grad}{\operatorname{grad}}

\title{Univalence Criteria for Lifts of Harmonic Mappings to Minimal Surfaces}

\thanks{The authors were supported in part by FONDECYT Grant \# 1030589.}

\author{M. Chuaqui \and P. Duren \and B. Osgood}

\address{P. Universidad Cat\'olica
de Chile}
\email{mchuaqui@mat.puc.cl}
\address{University of Michigan}
\email{duren@umich.edu}
\address{Stanford University}
\email{osgood@stanford.edu}

\subjclass[2000]{Primary 30C99; Secondary 31A05, 53A10}

\keywords{Harmonic mapping, Schwarzian derivative, curvature, 
minimal surface}

\date{}

\begin{document}

\maketitle

\bibliographystyle{amsplain}

\begin{abstract}A general criterion in terms of the Schwarzian derivative 
is given for global univalence of the Weierstrass--Enneper lift of a 
planar harmonic mapping.  Results on distortion and boundary regularity 
are also deduced. Examples are given to show that the criterion is sharp. The analysis depends on a generalized Schwarzian defined for conformal metrics and on a Schwarzian introduced by Ahlfors for curves. Convexity plays a central role.
\end{abstract}

\section{Introduction \label{section:intro}}

If a function $f$ is analytic and locally univalent, its 
{\em Schwarzian derivative} is defined by 
$$
\mathcal{S}f = (f''/f')' - \tfrac12(f''/f')^2 = f'''/f' 
- \tfrac32(f''/f')^2\,.
$$
The Schwarzian is invariant under postcomposition with M\"obius 
transformations
$$
T(z) = \frac{az+b}{cz+d}\,,\qquad ad-bc\neq 0\,;
$$
that is, $\mathcal{S}(T\circ f) = \mathcal{S}(f)$.  If $g$ is 
any function analytic and locally univalent on the range of $f$, then 
\begin{equation}
\mathcal{S}(g\circ f) = (\mathcal{S}(g)\circ f){f'}^2 + \mathcal{S}(f)\,.
\label{eq:chain-rule}
\end{equation}
In particular, $\mathcal{S}(g\circ T) = (\mathcal{S}(g)\circ T){T'}^2$, since 
$\mathcal{S}(T)=0$ for every M\"obius transformation $T$.  For an arbitrary 
analytic function $\psi$, the functions $f$ with Schwarzian 
$\mathcal{S}f=2\psi$ are those of the form $f=w_1/w_2$, where $w_1$ and $w_2$ 
are linearly independent solutions of the linear differential equation 
$w''+ \psi w=0$.  It follows that M\"obius transformations are the only 
functions with  $\mathcal{S}f=0$. More generally, if $\mathcal{S}f=\mathcal{S}g$, 
then $f=T\circ g$ for some M\"obius transformation $T$.

     In a groundbreaking paper, Nehari \cite{nehari:schlicht} showed that estimates 
on the Schwarzian provide criteria for global univalence.  He made the 
key observation that if a function $f$ is analytic and locally univalent 
in a simply connected domain $\Omega$, with Schwarzian $\mathcal{S}f=2\psi$, 
then $f$ is globally univalent if and only if no solution of the 
differential equation $w''+\psi w=0$, other than the zero solution, 
vanishes more than once in $\Omega$.  The univalence problem then reduces 
to a question about differential equations that can be analyzed by means of 
the Sturm comparison theorem.  Specifically, Nehari proved that if 
$f$ is analytic and locally univalent in the unit disk $\Bbb D$, and if 
\begin{equation}
|\mathcal{S}f(z)| \leq \frac{2}{(1-|z|^2)^2}\,,\qquad z\in \Bbb D\,, 
\label{eq:nehari-2}
\end{equation}
then $f$ is univalent in $\Bbb D$.  He also showed that the inequality 
\begin{equation}
|\mathcal{S}f(z)|\leq \frac{\pi^2}{2}
\label{eq:nehari-pi^2/2}
\end{equation}
 implies univalence.  Later he showed \cite{nehari:nehari-old-p} (see also \cite{nehari:nehari-p})  
that $f$ is univalent under the general hypothesis 
\begin{equation}
|\mathcal{S}f(z)|\leq 2p(|z|)\,,
\label{eq:classical-p}
\end{equation}
 where $p(x)$ is a positive, continuous, even 
function with the properties that $(1-x^2)^2p(x)$ is nonincreasing on 
the interval $[0,1)$ and no nontrivial solution $u$ of the differential 
equation $u''+pu=0$ has more than one zero in $(-1,1)$.  We will 
refer to such functions $p(x)$ as {\it Nehari functions}.  Conditions  \eqref{eq:nehari-2} and 
\eqref{eq:nehari-pi^2/2} are special cases of \eqref{eq:classical-p}.  Ess\'en and 
Keogh \cite{essen-keogh:schwarzian} solved some extremal problems for analytic functions 
satisfying Nehari's general condition.   Osgood and Stowe \cite{os:nehari} developed 
a common generalization of Nehari's criteria, and others, involving the 
curvature of a metric.  
     
     The main purpose of the present paper is to derive a corresponding 
univalence criterion for harmonic mappings, or rather for their lifts to 
minimal surfaces.  In previous work \cite{cdo:harmonic-schwarzian}, \cite{cdo:curvature}, \cite{cdo:ellipses} we have defined the 
Schwarzian derivative of a harmonic mapping and have developed some 
of its properties.  It is natural to identify a harmonic mapping with 
its Weierstrass--Enneper lift to a minimal surface, and it is this 
lift whose univalence is implied by bounds on the Schwarzian 
derivative.  For the underlying harmonic mappings, univalent or not, 
our criterion is also shown to imply estimates on distortion and 
properties of boundary regularity that are better than those known 
or conjectured (see \cite{sheil-small:harmonic-constants} or \cite{duren:harmonic}) for the full normalized class of 
univalent harmonic mappings.  In this respect our investigation can 
be viewed as a harmonic analogue of earlier work on analytic functions 
by Gehring and Pommerenke \cite{gehring:gehring-pommerenke} and Chuaqui and Osgood \cite{co:gp}, \cite{co:noncomplete}, \cite{co:extremal}.

     A planar harmonic mapping is a complex-valued harmonic function $f(z)$, $z=x+iy$, 
defined on some domain $\Omega\subset\Bbb C$.  If $\Omega$ is simply 
connected, the mapping has a canonical decomposition $f=h+\overline{g}$, 
where $h$ and $g$ are analytic in $\Omega$ and $g(z_0)=0$ for some specified 
point $z_0\in \Omega$.  The mapping $f$ is locally univalent if and only if 
its Jacobian $|h'|^2 - |g'|^2$ does not vanish.  It is said to be
orientation-preserving if $|h'(z)|>|g'(z)|$ in $\Omega$, or equivalently if
$h'(z)\neq0$ and the dilatation $\omega=g'/h'$ has the property
$|\omega(z)|<1$ in $\Omega$.

     According to the Weierstrass--Enneper formulas, a harmonic mapping 
$f=h+\overline{g}$ with $|h'(z)|+|g'(z)|\neq0$ lifts locally to map into 
a minimal surface, $\Sigma$, described by conformal parameters if and only if 
its dilatation $\omega=q^2$, the square of a meromorphic function $q$.  
The Cartesian coordinates $(U,V,W)$ of the surface are then given by
$$
U(z)=\text{Re}\{f(z)\}\,,\quad 
V(z)=\text{Im}\{f(z)\}\,,\quad
W(z)= 2\,\text{Im}\left\{\int_{z_0}^z h'(\zeta)q(\zeta)\,d\zeta\right\}\,.
$$
We use the notation 
$$
\widetilde{f}(z) = \bigl(U(z),V(z),W(z)\bigr)
$$
for the lifted mapping of $\Omega$ into $\Sigma$.  The height of the 
surface can be expressed more symmetrically as 
$$
W(z)= 2\,\text{Im}\left\{\int_{z_0}^z 
\sqrt{h'(\zeta)g'(\zeta)}\,d\zeta\right\}\,,
$$ 
since a requirement equivalent to $\omega=q^2$ is that $h'g'$ be the 
square of an analytic function.  The first fundamental form of the 
surface is $ds^2=e^{2\sigma}|dz|^2$, where the conformal factor is 
$$
e^\sigma = |h'|+|g'|\,.
$$
The Gauss curvature of the surface at a point $\widetilde{f}(z)$ for 
which $h'(z)\neq0$ is 
\begin{equation}
K = - e^{-2\sigma}\Delta \sigma
= - \frac{4|q'|^2}{|h'|^2(1+|q|^2)^4}\,,  
\label{eq:curvature}
\end{equation}  
where $\Delta$ is the Laplacian operator.  Further information about 
harmonic mappings and their relation to minimal surfaces can be found 
in \cite{duren:harmonic}.
 
   For a harmonic mapping $f=h+\overline{g}$ with $|h'(z)|+|g'(z)|\neq0$, 
whose dilatation is the square of a meromorphic function, we have 
defined \cite{cdo:harmonic-schwarzian} the {\em Schwarzian derivative} by the formula 
\begin{equation}
\mathcal{S}f = 2\bigl(\sigma_{zz} - \sigma_z^2\bigr)\,, 
\label{eq:harmonic-schwarzian}
\end{equation}
where 
$$
\sigma_z = \frac{\partial\sigma}{\partial z}  
= \frac12 \left(\frac{\partial\sigma}{\partial x} 
- i \frac{\partial\sigma}{\partial y}\right)\,, \qquad z = x+iy\,.
$$
Some background for this definition, on conformal Schwarzians, is discussed in Section \ref{sec:conformal-schwarzian}.
With $h'(z)\neq0$ and $g'/h'=q^2$, a calculation ({\it cf}. \cite{cdo:harmonic-schwarzian}) 
produces the expression
$$
\mathcal{S}f = \mathcal{S}h +\frac{2\overline{q}}{1+|q|^2}
\left(q'' - \frac{q'h''}{h'}\right) -4\left(\frac{q'\overline{q}}{1+|q|^2}
\right)^2\,.
$$ 
As observed in \cite{cdo:harmonic-schwarzian}, the formula remains valid if $\omega$ is 
not a perfect square, provided that neither $h'$ nor $g'$ has a simple 
zero.

     It must be emphasized that we are not requiring our harmonic 
mappings to be locally univalent.  In other words, the Jacobian need 
not be of constant sign in the domain $\Omega$.  The orientation of 
the mapping may reverse, corresponding to a folding in the associated 
minimal surface.  It is also possible for the minimal surface to 
exhibit several sheets above a point in the $(U,V)$--plane.  Thus the 
lifted mapping $\widetilde{f}$ may be univalent even when the underlying 
mapping $f$ is not.   

     The following theorem gives a criterion for the lift of a harmonic 
map to be univalent.

\begin{thm} \label{thm:p-criterion}  Let $f=h+\overline{g}$ be a harmonic mapping of the 
unit disk, with $e^{\sigma(z)}=|h'(z)|+|g'(z)|\neq0$ and dilatation 
$g'/h'=q^2$ for some meromorphic function $q$.  Let $\widetilde{f}$ denote 
the Weierstrass--Enneper lift of $f$ into a minimal surface $\Sigma$ with 
Gauss curvature $K=K(\widetilde{f}(z))$ at the point $\widetilde{f}(z)$.  
Suppose that 
\begin{equation}
|\mathcal{S}f(z)| + e^{2\sigma(z)} |K(\widetilde{f}(z))| \leq 2p(|z|)\,, 
\qquad z\in\Bbb D\,, 
\label{eq:p-criterion} 
\end{equation}
for some Nehari function $p$. Then $\widetilde{f}$ is univalent in $\Bbb D$.
\end{thm}

     When $f$ is analytic and locally univalent in $\Bbb D$, the result 
reduces to Nehari's theorem cited earlier, since the minimal surface 
$\Sigma$ is then a plane with $K=0$.

Theorem 1 is sharp, but to show this, and to 
describe the extremal mappings, we need to know that $\widetilde{f}$ 
has a continuous extension to the closed disk.  We state this now as 
a theorem, although in fact we will obtain stronger results on the 
smoothness of the boundary function.  

\begin{thm} \label{thm:boundary-extension} Under the hypotheses of Theorem 1, the 
Weierstrass--Enneper lift $\widetilde{f}$ has an extension to 
$\overline{\Bbb D}$ that is continuous with respect to the spherical 
metric.
\end{thm}

     According to Theorem 2, the condition \eqref{eq:p-criterion} implies that the 
extended mapping $\widetilde{f}$ sends the unit circle to a continuous 
closed curve $\Gamma$ lying on the surface $\Sigma$ in ${\Bbb R}^3\cup\{\infty\}$.  In fact, $\Gamma$ 
is a simple closed curve, or equivalently the mapping $\widetilde{f}$ 
is univalent in the closed  unit disk, except in special circumstances 
which we now describe.  A harmonic mapping $f$ satisfying \eqref{eq:p-criterion} is said to be  
{\it extremal}  if $\Gamma$ is {\em not} a simple closed curve; that is, if 
$\widetilde{f}(\zeta_1)=\widetilde{f}(\zeta_2)$ for some pair of 
distinct points $\zeta_1$ and $\zeta_2$ on the unit circle 
$\partial\Bbb D$.  One calls $P=f(\zeta_1) = f(\zeta_2)$ a {\em cut point} of $\Gamma$. The following theorem describes a characteristic 
property of extremal mappings.

\begin{thm} 
\label{thm:extremal}
Under the hypotheses of Theorem 1, suppose the 
closed curve $\Gamma = \widetilde{f}(\partial\Bbb D)$ is not simple and let $P$ be a cut point.  
Then there exists a Euclidean circle or line $C$ such that $C\setminus\{P\}$ is a line of curvature 
of $\widetilde{f}({\Bbb D})$ on the surface $\Sigma$. Furthermore, equality holds in \eqref{eq:p-criterion} along $\widetilde{f}^{-1}(C \setminus\{P\})$.     
\end{thm}

 Theorems 1 and 3 will be proved in Section 3, Theorem 2 in Section 6.    In the last section of the paper we will construct some examples 
of extremal mappings illustrating the phenomenon in Theorem 3 and showing that the criterion in Theorem 1 is sharp, in fact best possible in some particular cases. We will also show in that section that such extremal lifts must actually map the disk into either a catenoid or a plane; this will be a consequence of Theorem 3 and a purely differential geometric property of the catenoid. 

We think it is striking that the theorems for analytic functions generalize in this manner, and one cannot help but wonder what other related aspects of classical geometric function theory have counterparts for harmonic mappings or their lifts. Of additional interest is that our analysis involves not only the classical Schwarzian and its conformal generalization, but also a version of the Schwarzian introduced by Ahlfors for curves in ${\Bbb R}^n$, to which we now turn.

\section{Ahlfors' Schwarzian and univalence along curves}

     Ahlfors \cite{ahlfors:schwarzian-rn} introduced a notion of Schwarzian derivative 
for mappings of a real interval into ${\Bbb R}^n$ by formulating 
suitable analogues of the real and imaginary parts of $\mathcal{S}f$ for 
analytic functions $f$.  A simple calculation shows that 
$$
\text{Re}\{\mathcal{S}f\} = \frac{\text{Re}\{f'''\overline{f'}\}}{|f'|^2} 
-3 \,\frac{\text{Re}\{f''\overline{f'}\}^2}{|f'|^4} + \frac32 
\frac{|f''|^2}{|f'|^2}\,.
$$
For mappings $\varphi : (a,b) \rightarrow {\Bbb R}^n$ of class $C^3$ with 
$\varphi'(x)\neq0$, Ahlfors defined the analogous expression
\begin{equation}
\label{eq:S1}
S_1\varphi = \frac{\langle \varphi''',\varphi'\rangle}{|\varphi'|^2}
- 3\frac{\langle \varphi'',\varphi'\rangle^2}{|\varphi'|^4} 
+ \frac32\frac{|\varphi''|^2}{|\varphi'|^2}\,, 
\end{equation}
where $\langle\cdot,\cdot\rangle$ denotes the Euclidean inner product 
and now $|{\bold x}|^2=\langle{\bold x},{\bold x}\rangle$ for ${\bold x}\in
{\Bbb R}^n$.  Ahlfors also defined a second expression analogous to 
$\text{Im}\{\mathcal{S}f\}$, but this is not relevant to the present 
discussion.  

     Ahlfors' Schwarzian is invariant under postcomposition with M\"obius 
transformations; that is, under every composition of rotations, 
magnifications, translations, and inversions in ${\Bbb R}^n$.  Only its 
invariance under inversion
$$
{\bold x} \mapsto \frac{\bold x}{|\bold x|^2}\,, \qquad 
{\bold x}\in{\Bbb R}^n\,,
$$
presents a difficulty; this can be checked by straightforward but tedious 
calculation.  It should also be noted that $S_1$ transforms as expected  
under change of parameters.  If $x=x(t)$ is a smooth function with 
$x'(t)\neq0$, and $\psi(t)=\varphi(x(t))$, then
$$
S_1\psi(t) = S_1\varphi(x(t))\,x'(t)^2 + \mathcal{S}x(t)\,.
$$ 

    With the notation $v=|\varphi'|$, Chuaqui and Gevirtz \cite{chuaqui-gevirtz:S1} used the Frenet-Serret formulas to show
that 
\begin{equation}
S_1\varphi = (v'/v)' - \tfrac12 (v'/v)^2 
+ \tfrac12 v^2 \kappa^2 = \mathcal{S}(s) + \tfrac12 v^2 \kappa^2\,, 
\label{eq:S1-and-S}
\end{equation}
where $s=s(x)$ is the arclength of the curve and $\kappa$ is its curvature.  Our proof of Theorem 1 will be based on    the following result, also due to Chuaqui and Gevirtz in \cite{chuaqui-gevirtz:S1}. 
\begin{thma} Let $p(x)$ be a continuous function such that the  
differential equation $u''(x)+p(x)u(x)=0$ admits no nontrivial solution 
$u(x)$ with more than one zero in $(-1,1)$. Let $\varphi : (-1,1)\rightarrow 
{\Bbb R}^n$ be a curve of class $C^3$ with tangent vector 
$\varphi'(x)\neq 0$.  If $S_1\varphi(x)\leq2p(x)$, then $\varphi$ is 
univalent.   
\end{thma}

     If the function $p(x)$ of Theorem A is even, as will be the case for a Nehari function, then the solution $u_0$ 
of the differential equation $u''+pu=0$ with initial conditions $u_0(0)=1$ 
and $u_0'(0)=0$ is also even, and therefore $u_0(x)\neq0$ on $(-1,1)$, 
since otherwise it would have at least two zeros.  Thus the function 
\begin{equation}
\Phi(x) = \int_0^x u_0(t)^{-2}\,dt\,, \qquad -1 < x < 1\,,
\label{eq:Phi}
\end{equation}
is well defined and has the properties $\Phi(0)=0$, $\Phi'(0)=1$, 
$\Phi''(0)=0$, $\Phi(-x)=-\Phi(x)$.  The standard method of reduction 
of order produces the independent solution $u=u_0\Phi$ to $u''+pu=0$, 
and so $\mathcal{S}\Phi=2p$. Note also that $S_1\Phi=\mathcal{S}\Phi$, since 
$\Phi$ is real-valued.  Thus $S_1\Phi=2p$.  

     The next theorem, again to be found in  \cite{chuaqui-gevirtz:S1}, asserts that 
the mapping $\Phi: (-1,1)\rightarrow {\Bbb R}\subset{\Bbb R}^n$ is extremal 
for Theorem A if $\Phi(1)=\infty$, and that every extremal mapping 
$\varphi$ is then a M\"obius postcomposition of $\Phi$.  

\begin{thmb} Let $p(x)$ be an even function with the properties 
assumed in Theorem A, and let $\Phi$ be defined as above.  Let 
$\varphi: (-1,1)\rightarrow {\Bbb R}^n$ satisfy $S_1\varphi(x)\leq2p(x)$ 
and have the normalization $\varphi(0)=0$, $|\varphi'(0)|=1$, and 
$\varphi''(0)=0$.  Then $|\varphi'(x)|\leq\Phi'(|x|)$ for $x\in(-1,1)$, 
and $\varphi$ has an extension to the closed interval $[-1,1]$ that is 
continuous with respect to the spherical metric.  Furthermore, there are 
two possibilities, as follows. 
\par $(i)$  If $\Phi(1)<\infty$, then $\varphi$ is univalent in $[-1,1]$ 
and $\varphi([-1,1])$ has finite length.
\par $(ii)$ If $\Phi(1)=\infty$, then either $\varphi$ is univalent in 
$[-1,1]$ or $\varphi=R\circ\Phi$ for some rotation $R$ of ${\Bbb R}^n$.  
\end{thmb}

     Note that in case $(ii)$ the mapping $\Phi$ sends both ends of the 
interval to the point at infinity and is therefore not univalent in 
$[-1,1]$. The role of $\Phi$ as an extremal for the harmonic univalence criterion \eqref{eq:p-criterion} will emerge in the following sections.

\section{Univalence and extremal functions: Proofs of Theorems 1 and 3}

\label{section:univalence}

   The proofs of Theorems 1 and 3 will be based on Theorems A and B.  
The following lemma makes the connection.

\begin{lem} Let $f$ be a harmonic mapping of the unit disk with 
nonvanishing conformal factor $e^\sigma$ and dilatation $q^2$ for some 
meromorphic function $q$.  Let $\widetilde{f}$ be the lift of $f$ to 
a minimal surface $\Sigma$ with Gauss curvature $K$.  Then 
$$
\begin{aligned}
S_1\widetilde{f}(x) &= {\rm Re}\{{\mathcal{S}}f(x)\} 
+ \tfrac12 e^{2\sigma(x)}|K(\widetilde{f}(x))| 
+ \tfrac12 e^{2\sigma(x)} \kappa_e(\widetilde{f}(x))^2 \\
&\leq {\rm Re}\{{\mathcal{S}}f(x)\} 
+ e^{2\sigma(x)}|K(\widetilde{f}(x))|\,, \qquad -1<x<1\,, 
\end{aligned}
$$
where $\kappa_e(\widetilde{f}(x))$ denotes the normal curvature of the 
curve $\widetilde{f}:(-1,1)\rightarrow \Sigma$ at the point $\widetilde{f}(x)$.  
Equality occurs at a point $x$ if and only if the curve is tangent to a 
line of curvature of \,$\Sigma$ at the point $\widetilde{f}(x)$.  
\end{lem}

\begin{proof} According to the formula \eqref{eq:S1-and-S}, we have
$$
S_1\widetilde{f} = (v'/v)' - \tfrac12 (v'/v)^2 + \tfrac12 v^2 \kappa^2\,,
$$
where $v=e^\sigma$ is the conformal factor of the surface $\Sigma$ and
$\kappa$ is the curvature of the curve $\widetilde{f}:(-1,1) \rightarrow
\Sigma$.
The tangential and normal projections of the curvature vector
$$
\frac{d{\bold t}}{ds} = \frac{d}{ds}\left(\frac{d\widetilde{f}}{ds}\right)
$$
are the geodesic or intrinsic curvature $\kappa_i$ and the normal or 
extrinsic curvature $\kappa_e$, respectively.  Thus 
$\kappa^2 = \kappa_i^2 + \kappa_e^2$.  In a previous paper \cite{cdo:curvature} we related 
the geodesic curvature of the lifted curve to the curvature of an 
underlying curve $C$ in the parametric plane.  In the present context $C$ 
is the linear segment $(-1,1)$, with curvature zero, and so our formula 
(\cite{cdo:curvature}, eq. (4)) reduces to $e^\sigma \kappa_i = -\sigma_y$, where
$\sigma_y = \partial\sigma/\partial y$.  Therefore, we find that
$$
\begin{aligned}
S_1\widetilde{f} &= \sigma_{xx} - \tfrac12 \sigma_x^2 + \tfrac12
e^{2\sigma}
\bigl(\kappa_i^2 + \kappa_e^2\bigr) \\
&= \tfrac12\bigl(\sigma_{xx} - \sigma_{yy} - \sigma_x^2 + \sigma_y^2\bigr)
+ \tfrac12\bigl(\sigma_{xx} + \sigma_{yy}\bigr) + \tfrac12 e^{2\sigma}
\kappa_e^2 \\
&= \text{Re}\{{\mathcal{S}}f\} + \tfrac12 e^{2\sigma}|K| + \tfrac12 e^{2\sigma}
\kappa_e^2 \,,
\end{aligned}
$$
in view of the expression \eqref{eq:curvature} for the Gauss curvature $K$.  Here we have
used the fact that $K\leq0$ for a minimal surface.  In fact, $K=k_1k_2$,
where $k_1$ and $k_2$ are the principal curvatures, the maximum and
minimum
of the normal curvature $\kappa_e$ of the surface as the tangent direction
varies.  But a minimal surface has mean curvature zero; that is, 
$k_1 + k_2 = 0$.  This implies that $K\leq0$ and 
$\kappa_e^2\leq|k_1k_2|=|K|$.  Consequently,
$$
S_1\widetilde{f} \leq \text{Re}\{{\mathcal{S}}f\} + e^{2\sigma}|K|\,,
$$
with equality if and only if the curve is tangent to a line of curvature
of the surface $\Sigma$.  In other words, equality occurs precisely when
the curve is aligned in a direction of maximum or minimum normal
curvature.
\end{proof}

\begin{proof}[Proof of Theorem 1]  By Lemma 1, the inequality \eqref{eq:p-criterion} implies that 
the curve $\widetilde{f}: (-1,1)\rightarrow \Sigma$ satisfies the hypothesis 
$S_1\widetilde{f}(x) \leq 2p(x)$ of Theorem A.  Thus Theorem A tells us that 
$\widetilde{f}$ is univalent in the interval $(-1,1)$.  

     In order to conclude that $\widetilde{f}$ is univalent in $\Bbb D$, 
we adapt a clever argument due to Nehari \cite{nehari:nehari-old-p}.  We want to show that 
$\widetilde{f}(z_1)\neq\widetilde{f}(z_2)$ for any given pair of distinct 
points $z_1,z_2\in\Bbb D$.  But if $f$ satisfies \eqref{eq:p-criterion}, then so does every 
rotation $f(e^{i\theta}z)$.  Consequently, we may assume that the 
hyperbolic geodesic passing through the points $z_1$ and $z_2$ intersects 
the imaginary axis orthogonally.  Let $i\rho$ denote the point of 
intersection, and observe that the M\"obius transformation 
\begin{equation}
T(z) = \frac{i\rho - z}{1+i\rho z}\,, \qquad z\in\Bbb D\,,  
\label{eq:special-mobius}
\end{equation}
maps the disk onto itself and preserves the imaginary axis, so that it 
maps the given geodesic onto the real segment $(-1,1)$.  Moreover, $T$ 
is an involution, so that $T=T^{-1}$ and $T$ also maps the segment 
$(-1,1)$ onto the geodesic through $z_1$ and $z_2$.  In particular, 
$T(x_1)=z_1$ and $T(x_2)=z_2$ for some points $x_1$ and $x_2$ in the 
interval $(-1,1)$.  The composite function $F(z)=f(T(z))$ is a harmonic 
mapping of the disk whose lift $\widetilde{F}$ again maps $\Bbb D$ to 
the minimal surface $\Sigma$.  We claim that 
\begin{equation}
|{\mathcal{S}}F(x)| + e^{2\tau(x)} |K(\widetilde{F}(x))| \leq 2p(x)\,, 
\qquad -1 < x< 1\,, 
\label{eq:SF}
\end{equation}
where $e^\tau=|H'|+|G'|$ is the conformal factor associated with 
$F=H+\overline{G}$.  Indeed, by the chain rule,
$$
e^{\tau(x)} = e^{\sigma(T(x))}|T'(x)|\,,
$$
whereas ${\mathcal{S}}F(x) = {\mathcal{S}}f(T(x))\, T'(x)^2$ and 
$\widetilde{F}(x)=\widetilde{f}(T(x))$.  Therefore, by virtue of the 
hypothesis \eqref{eq:p-criterion}, the claim \eqref{eq:SF} will be established if we can show that 
\begin{equation}
|T'(x)|^2 p(|T(x)|) \leq p(x)\,, \qquad -1 < x< 1\,. 
\label{eq:T'<=p}
\end{equation}
But 
$$
\frac{|T'(z)|}{1 - |T(z)|^2} = \frac{1}{1-|z|^2}\,, \qquad z\in\Bbb D\,,
$$
and so \eqref{eq:T'<=p} reduces to the inequality 
\begin{equation}
\left(1 - |T(x)|^2\right)^2  p(|T(x)|) \leq \left(1-x^2\right)^2 p(x)\,,
\qquad -1<x<1\,, 
\label{eq:nehari-trick}
\end{equation}
which follows from the assumption that $\left(1-x^2\right)^2 p(x)$ is 
nonincreasing on $[0,1)$, since $p$ is an even function and an easy 
calculation shows that $|T(x)|>|x|$.  This proves our claim \eqref{eq:SF}.  

     Finally, we return to the remark made at the beginning of the proof.  
In view of Lemma 1 and Theorem A, the inequality \eqref{eq:SF} implies that 
$\widetilde{F}$ is univalent in $(-1,1)$.  Therefore, 
$\widetilde{F}(x_1)\neq\widetilde{F}(x_2)$, which says that 
$\widetilde{f}(z_1)\neq\widetilde{f}(z_2)$.  This proves Theorem 1.
\end{proof}

\begin{proof}[Proof of Theorem 3]  We start with the assumption that a harmonic 
mapping $f$ satisfies the hypotheses of Theorem 1 for some Nehari function 
$p(x)$, and its lift $\widetilde{f}$ has the property 
$\widetilde{f}(\zeta_1)=\widetilde{f}(\zeta_2)$ for some pair of distinct 
points $\zeta_1$ and $\zeta_2$ on the unit circle.  The first step  
is to show that after suitable modification of $f$ we may 
take $\zeta_1$ and $\zeta_2$ to be diametrically opposite points.  
More precisely, we will show that some M\"obius transformation $T$ of 
$\Bbb D$ onto itself produces a harmonic mapping $F=f\circ T$ with a 
lift $\widetilde{F}(z)= \widetilde{f}(T(z))$ for which the 
inequality \eqref{eq:SF} holds and $\widetilde{F}(1)=\widetilde{F}(-1)$.  To make 
the analysis clear, we will distinguish two cases.

     Suppose first that $(1-x^2)^2p(x)$ is constant.  Then equality 
holds in \eqref{eq:nehari-trick}, hence in \eqref{eq:T'<=p}, for every M\"obius automorphism of 
$\Bbb D$.  Consequently, every mapping $F=f\circ T$ satisfies the 
inequality \eqref{eq:SF}, and we will have $\widetilde{F}(1)=\widetilde{F}(-1)$ 
if we choose the automorphism $T$ such that $T(1)=\zeta_1$ and 
$T(-1)=\zeta_2$.  

     If $(1-x^2)^2p(x)$ is not constant, we claim that necessarily 
$\zeta_1=-\zeta_2$.  If not, then after suitable rotation we may 
assume that 
$$
\text{Im}\{\zeta_1\}=\text{Im}\{\zeta_2\}>0 \qquad \text{and} \qquad 
\text{Re}\{\zeta_2\}=-\text{Re}\{\zeta_1\} >0\,.
$$
Now let $i\rho$ be the point where the imaginary axis meets the 
hyperbolic geodesic from $\zeta_1$ to $\zeta_2$, and observe that the 
M\"obius transformation \eqref{eq:special-mobius} carries the real segment $(-1,1)$ onto 
this geodesic, with $T(1)=\zeta_1$ and  $T(-1)=\zeta_2$.  As in the 
proof of Theorem 1, we arrive at the inequality \eqref{eq:nehari-trick}, which implies 
\eqref{eq:T'<=p}  and therefore \eqref{eq:SF}, but now the inequality \eqref{eq:SF} cannot reduce to 
equality throughout the interval $(-1,1)$, since $(1-x^2)^2p(x)$ is not 
constant.  In view of Lemma 1, we conclude that 
$$
S_1\widetilde{F}(x) \leq 2 p(x)\,, \qquad -1<x<1\,,
$$
with strict inequality in some part of the interval.  In particular, 
$S_1\widetilde{F}\neq 2p$.  However, this conclusion stands in 
contradiction to Theorem B.  Indeed, after postcomposition with a 
suitable M\"obius transformation $M$ of ${\Bbb R}^3$, we obtain a 
mapping 
$$
\varphi = M\circ\widetilde{F} : (-1,1) \rightarrow{\Bbb R}^3 
$$
with the required normalization $\varphi(0)=0$, $|\varphi'(0)|=1$, and 
$\varphi''(0)=0$ ({\it cf.} \cite{chuaqui-gevirtz:S1}).  But $\varphi(1)=\varphi(-1)$, so 
Theorem B tells us that $\varphi = R\circ\Phi$ for some rotation $R$.  
Hence $(R^{-1}M)\circ\widetilde{F}=\Phi$, and so $S_1\widetilde{F}
=S_1\Phi=2p$.  This contradicts our earlier conclusion and shows that 
$\zeta_1=-\zeta_2$.  Thus some rotation $T$ of the disk produces a 
harmonic mapping $F=f\circ T$ whose lift 
$\widetilde{F}=\widetilde{f}\circ T$ satisfies \eqref{eq:SF} and 
$\widetilde{F}(1)=\widetilde{F}(-1)$.

     In all cases we find that some M\"obius transformation $T$ of the 
disk onto itself produces a harmonic mapping $F=f\circ T$ that satisfies 
the inequality \eqref{eq:SF}, and whose lift $\widetilde{F}=\widetilde{f}\circ T$ 
has the property $\widetilde{F}(1)=\widetilde{F}(-1)$.  Thus 
$\widetilde{F}$ maps the interval $[-1,1]$ to a closed curve on the 
surface $\Sigma$.  As previously indicated, Theorem B then shows that 
$\widetilde{F}=V\circ\Phi$ for some M\"obius transformation $V$ of 
${\Bbb R}^3$.  It follows that $S_1\widetilde{F}=S_1\Phi=2p$, and also 
that $\widetilde{F}$ maps $[-1,1]$ to a Euclidean circle or line, since 
$\Phi$ maps $[-1,1]$ onto the extended real line and M\"obius 
transformations preserve circles.  On the other hand, since $F$ 
satisfies \eqref{eq:SF}, Lemma 1 shows that 
$$
S_1\widetilde{F}(x) \leq \text{Re}\{{\mathcal{S}}F(x)\} 
+ e^{2\tau(x)}|K(\widetilde{F}(x))| \leq 2p(x)\,, \qquad -1<x<1\,.
$$
But $S_1\widetilde{F}(x)=2p(x)$, so equality holds throughout.  
According to Lemma 1, this says that the circle 
$\widetilde{F}([-1,1])$ is everywhere tangent to a line of curvature, 
so it is in fact a line of curvature of $\Sigma$.
\end{proof}

\section{Conformal Schwarzian}
 \label{sec:conformal-schwarzian}

Results on extensions to the boundary and estimates on distortion for harmonic mappings satisfying the univalence criterion \eqref{eq:p-criterion} depend upon inequalities derived from convexity. These in turn call on a generalized Schwarzian that is computed with respect to a conformal metric and on a second order differential equation associated with the Schwarzian. It is this `conformal Schwarzian' when specialized  to the lift of a harmonic mapping that produces the definition 
\eqref{eq:harmonic-schwarzian}; see also \cite{cdo:harmonic-schwarzian}. The definition and properties are suggested by the classical case, and have analogues there, but the generalization must be framed in the terminology of differential geometry; see for example \cite{oneill:semiriemannian} for a very accessible treatment of the operations we use here. This section provides a brief summary of the generalized Schwarzian, with all definitions given for dimension 2.  We refer to \cite{os:sch} for the higher dimensional setting, and to \cite{co:extremal} for applications of convexity in 2 dimensions, similar to what we will do here for harmonic mappings.

Let $\mathbf{g}$ be a Riemannian metric on the disk $\Bbb D$. We may assume that $\mathbf{g}$ is conformal to the Euclidean metric, $\mathbf{g}_0=dx\otimes dx +dy \otimes dy= |dz|^2$. Let $\psi$ be a smooth function on $\Bbb D$ and form the symmetric 2-tensor
\begin{equation}
\Hess_{\mathbf{g}}(\psi) - d\psi \otimes d\psi.
\label{eq:Hessian}
\end{equation}
Here $\Hess$ denotes the Hessian operator. For example, if $\gamma(s)$ is an arclength parametrized geodesic for $\mathbf{g}$, then
\[
\Hess_{\mathbf{g}}(\psi)(\gamma',\gamma') = \frac{d^2}{ds^2}(\psi\circ \gamma)\,.
\]
The Hessian depends on the metric, and since we will be changing metrics we indicate this dependence by the subscript $\mathbf{g}$. 

With some imagination the tensor \eqref{eq:Hessian} begins to resemble a Schwarzian; among other occurrences in differential geometry, it arises (in  2 dimensions) if one differentiates the equation that relates the geodesic curvatures of a curve for two conformal metrics. Such a curvature formula is a classical interpretation of the Schwarzian derivative, see \cite{os:sch} and \cite{cdo:curvature}. The trace of the tensor is the function
\[
\frac{1}{2}(\Delta_\mathbf{g} \psi - ||\grad_\mathbf{g}\psi||_\mathbf{g}^2),
\]
where again we have indicated by a subscript that the Laplacian, gradient and norm all depend on $\mathbf{g}$. It turns out to be most convenient to work with a traceless tensor when generalizing the Schwarzian, so we subtract off this function times the metric $\mathbf{g}$ and define the \textit{Schwarzian tensor} to be the symmetric, traceless, 2-tensor
\[ 
B_\mathbf{g}(\psi)={\Hess}_\mathbf{g}(\psi)-d\psi\otimes d\psi-\frac{1}{2}(\Delta_\mathbf{g} \psi-||\grad_\mathbf{g} \psi||^2)\mathbf{g}\,.
\label{eq:B-psi}
\]
Working in standard Cartesian coordinates one can represent  
$B_{\mathbf{g}}(\psi)$ as a symmetric, traceless $2\times 2$ matrix, say of the form 
\[
\begin{pmatrix}
a & -b \\
-b & -a
\end{pmatrix}\,.
\]
Further identifying such a matrix with the complex number $a+bi$ then allows us to associate the tensor $B_{\mathbf{g}}(\psi)$ with $a+bi$.

At each point $z\in\Bbb D$, the expression $B_\mathbf{g}(\psi)(z)$ is a bilinear form on the tangent space at $z$, and so its norm is
\[
||B_\mathbf{g}(\psi)(z)||_\mathbf{g} =\sup_{X,Y}B_\mathbf{g}(\psi)(z)(X,Y)\,,
\]
where the supremum is over unit vectors in the metric $\mathbf{g}$. If we compute the tensor with respect to the Euclidean metric and make the identification with a complex number as above, then
\[
||B_\mathbf{g_0}(\psi)(z)||_\mathbf{g_0} = |a+bi|\,.
\]

Now, if $f$ is analytic and locally univalent in $\Bbb D$, then it is a conformal mapping of $\Bbb D$ with the metric $\mathbf{g}$ into $\Bbb C$ with the Euclidean metric. The pullback $f^*\mathbf{g}_0$ is a metric on ${\Bbb D}$ conformal to $\mathbf{g}$, say $f^*\mathbf{g}_0 = e^{2\psi}\mathbf{g}$, and  the (conformal) Schwarzian of $f$ is now defined to be
\[
\mathcal{S}_\mathbf{g} f = B_\mathbf{g}(\psi)\,.
\]
If we take $\mathbf{g}$ to be the Euclidean metric then $\psi = \log|f'|$. Computing $B_{\mathbf{g}_0}(\log|f'|)$ and writing it in matrix form as above results in
\[
B_{\mathbf{g}_0}(\log|f'|)=\left( \begin{array}{rr} {\rm Re}\,\mathcal{S}f & -{\rm Im}\,\mathcal{S}f \\ -{\rm Im}\,\mathcal{S}f & -{\rm Re}\,\mathcal{S}f
\end{array}\right) \, ,
\]
where $\mathcal{S}f$ is the classical Schwarzian derivative of $f$. In this way we identify $B_{\mathbf{g}_0}(\log|f'|)$ with  $\mathcal{S}f$. 
 
Next, if $f=h+\overline{g}$ is a harmonic mapping of $\Bbb D$ and $\sigma = \log(|h'|+|g'|)$ is the conformal factor associated with the lift $\widetilde{f}$, we put
\[
\mathcal{S}f = \mathcal{S}_{\mathbf{g}_0} \widetilde{f} = B_{\mathbf{g}_0}(\sigma).
\]
Calculating this out and making the identification of the generalized Schwarzian with a complex number produces
\[
B_{\mathbf{g}_0}(\sigma) = 2(\sigma_{zz} - \sigma_z^2)\,,
\]
which is the definition of $\mathcal{S}f$ given in \eqref{eq:harmonic-schwarzian}. 

We need two more general facts. First, if we change a metric $\mathbf{g}$ conformally to $\widehat{\mathbf{g}} = e^{2\rho}\mathbf{g}$ then the tensor $B_\mathbf{g}(\psi)$ changes according to
\[
B_\mathbf{g}(\rho+\psi) = B_\mathbf{g}(\rho)+B_{\widehat{\mathbf{g}}}(\psi).
\]
This is actually a generalization of the chain rule \eqref{eq:chain-rule} for the Schwarzian.  An equivalent formulation is
\begin{equation}
B_{\widehat{\mathbf{g}}}(\psi -\rho) = B_\mathbf{g}(\psi) - B_\mathbf{g}(\rho)\,,
\label{eq:subtraction-formula}
\end{equation}
which is what we will need in the next section.

\medskip

Second, just as the linear differential equation $w''+p w=0$ is associated with $Sf=2p$, so is there a linear differential equation associated with the Schwarzian tensor. If $B_\mathbf{g}(\psi) = p$, where $p$ is a symmetric, traceless 2-tensor, then $\phi=e^{-\psi}$$\phi$ satisfies
\begin{equation}
\Hess_\mathbf{g}(\phi) + \phi p = \frac{1}{2}(\Delta_\mathbf{g}\phi)\mathbf{g}\,.
\label{eq:Hess-eq}
\end{equation}
Although the setting is more general, the substitution $\phi=e^{-\psi}$ is suggested by the classical case; regard $\mathcal{S}f=2p$ as a Riccati equation $u'-(1/2)u^2 = 2p$ for $u=f''/f'$. With $v=\log f'$ and $w=e^{-v/2} = (f')^{-1/2}$ one finds that $w$ satisfies $w''+pw=0$. See \cite{osgood:old-and-new}.

\medskip

Finally, a comment about convexity. Let $\gamma(s)$  be an arc-length parametrized geodesic for the metric $\mathbf{g}$. Evaluating both sides of the equation \eqref{eq:Hess-eq}
at the pair $(\gamma'(s),\gamma'(s))$ gives the scalar equation
\[
\frac{d^2}{ds^2}\phi(\gamma(s)) +\phi(\gamma(s)) p(\gamma'(s),\gamma'(s)) =\frac{1}{2} \Delta_\mathbf{g} \phi(\gamma(s))\,;
\]
this uses $\mathbf{g}(\gamma'(s),\gamma'(s))=1$. If  
\[
\frac{1}{2} \Delta_\mathbf{g} \phi(\gamma(s)) - \phi(\gamma(s)) p(\gamma'(s),\gamma'(s)) \ge 0
\]
for all geodesics then $\phi$ is convex relative to the metric $\mathbf{g}$. Without evaluating on a pair of vector fields, the condition for a function $\phi$ to be convex can be written as
\[
\Hess_\mathbf{g} \phi \ge \alpha\mathbf{g}\,,
\]
where $\alpha$ is a nonnegative function.  We will find that an upper bound for 
$\mathcal{S}_\mathbf{g}f$ coming from the univalence criterion \eqref{eq:p-criterion} leads via \eqref{eq:subtraction-formula} and \eqref{eq:Hess-eq} to just such a positive lower bound for the Hessian of an associated function.

\section{Univalence criteria and convexity} \label{section:convexity}

Convexity enters the picture by relating upper bounds on the Schwarzian tensor, in the guise of the Schwarzian of a harmonic mapping, to lower bounds on the Hessian of an associated function. There are two aspects to this. The first is to identify the appropriate  conformal metric with respect to which the computations are made, and we do this now.

We recall the function 
\[
\Phi(x) = \int_0^x u_0(t)^{-2}\,dt\,, \qquad -1 < x < 1\,,
\]
defined in \eqref{eq:Phi}, where $u_0$ is a positive solution of
\[
u''+ pu =0, \quad u(0)=1, u'(0)=0\,,
\]
when $p$ is a Nehari function. Recall also that $\Phi$ is odd with $\Phi(0)=0$, $\Phi'(0)=1$ and $\Phi''(0)=0$.

We use $\Phi$ to form the radial conformal metric $\mathbf{g} = \Phi'(|z|)^2|dz|^2$ on $\Bbb D$. It is straightforward to express the curvature as  
\begin{equation}
K_\mathbf{g}(z) = -2|\Phi'(r)|^{-2}(A(r)+p(r))\,,\qquad r=|z|\,,
\end{equation}
where
\begin{equation}
A(r) = \frac{1}{4}\left(\frac{\Phi''(r)}{\Phi'(r)}\right)^2 +\frac{1}{2r}\frac{\Phi''(r)}{\Phi'(r)}\,;
\label{eq:A(r)}
\end{equation}
see \cite{co:noncomplete}. Note also that $A(r)$ is continuous at $0$ with $A(0)=p(0)$, and that the curvature is {\em negative}. In particular
\begin{equation}
|K_\mathbf{g}(z)| = 2|\Phi'(r)|^{-2}(A(r)+p(r))\,.
\label{eq:|K_g|}
\end{equation}

Appealing to the results in \cite{co:noncomplete} we can assume that the metric $\Phi'(|z|)|dz|$ is complete, which means precisely that $\Phi(1) = \infty$. To elaborate, if $\Phi(1)<\infty$ then, as is shown in \cite{co:noncomplete}, there is a maximum value $t_0>1$ such that $t_0p(x)$ remains a Nehari function and such that the corresponding extremal $\Phi_{t_0}$ has $\Phi_{t_0}(1)=\infty$. Since any condition of the form $|\mathcal{S}f|+\cdots \le 2p$ implies trivially that $|\mathcal{S}f|+\cdots \le 2t_0p$ one may take $\Phi(1)=\infty$ at the outset. We make this assumption.  

Geometric consequences of completeness of the metric are that any two points in ${\Bbb D}$ can be joined by a geodesic for $\mathbf{g}$, and that any geodesic can be extended indefinitely. An analytic consequence of completeness, following  from $\Phi(1)=\infty$, is 
\begin{equation}
p(x)\le A(x)\,.
\label{eq:p<=A}
\end{equation}
This was shown in \cite{co:noncomplete}.

In Theorem 3 we have already seen $\Phi$ play the role of an extremal function for the univalence criterion. Our results on distortion and boundary behavior depend on $\Phi$ defining, as above, an `extremal metric' for the criterion. In the planar, analytic case a detailed study was carried out in \cite{co:noncomplete}  and \cite{co:extremal}.  

\medskip

The second aspect of our analysis is captured in the following convexity result.

\begin{thm} \label{prop:convexity} Let $f=h+\overline{g}$ be a harmonic
mapping satisfying the univalence criterion \eqref{eq:p-criterion},
\[
|\mathcal{S}f(z)|+e^{2\sigma(z)}|K(\widetilde{f}(z))| \leq 2p(|z|) = S\Phi(|z|) \,. 
\label{eq:criterion-for-convexity} 
\]
Let $\vphi(z) = \log\Phi'(|z|)$ and define
\begin{equation}
u(z)= e^{(\vphi(z)-\sigma(z))/2} = \sqrt{\frac{\Phi'(|z|)}{|h'(z)|+|g'(z)|}} \,. \label{eq:u-for-convexity}
\end{equation}
Then
\begin{equation}
{\Hess}_{\mathbf{g}}(u) \geq \frac{1}{4}u^{-3}|K|\mathbf{g} \, , \label{eq:bound-on-hessian}
\end{equation}
where $\mathbf{g}$ is the conformal metric $\Phi'(|z|)^2|dz|^2$. In
particular, $u$ is a convex function relative to $\mathbf{g}$.
\end{thm}

Thus we see that convexity obtains not for the conformal factor of $\widetilde{f}$, or for its square root, but for the square root of the {\em ratio} of the conformal factors of $\widetilde{f}$ and the extremal mapping for the univalence criterion. 

The work of the present section is to establish Theorem 4. The inequality \eqref{eq:p<=A} is crucial, and to highlight this aspect of the proof we split off the following calculation as a separate lemma. 

\begin{lem} 
If $f$ satisfies \eqref{eq:p-criterion} then
\begin{equation}
||B_{\mathbf{g}}(\sigma-\varphi)(z)||_{\mathbf{g}}+ e^{2(\sigma-\varphi)(z)}|K(\widetilde{f}(z))| \leq
\frac{1}{2}|K_\mathbf{g}(z)|\,. \label{eq:Bg(sigma-phi)}
\end{equation}
\end{lem}

\begin{proof} Using the variant \eqref{eq:subtraction-formula} of the addition formula for the Schwarzian tensor we have
\[
B_\mathbf{g}(\sigma-\vphi) = B_{\mathbf{g}_0}(\sigma) - B_{\mathbf{g}_0}(\vphi)\,,
\]
where $\mathbf{g}_0$ is the Euclidean metric. Now $\mathbf{g} = e^{2\vphi}\mathbf{g_0}$, and so the norm scales according to
\[
e^{2\vphi(z)}||B_\mathbf{g}(\sigma-\vphi)(z)||_\mathbf{g} = ||B_\mathbf{g}(\sigma-\vphi)(z)||_\mathbf{g_0}\,,
\]
whence by the preceding equation
\[
e^{2\vphi(z)}||B_\mathbf{g}(\sigma-\vphi)(z)||_\mathbf{g} =||B_{\mathbf{g}_0}(\sigma)(z) - B_{\mathbf{g}_0}(\vphi)(z)||_{\mathbf{g}_0}\,.
\]
Finally, 
\[
||B_{\mathbf{g}_0}(\sigma)(z) - B_{\mathbf{g}_0}(\vphi)(z)||_{\mathbf{g}_0} =|B_{\mathbf{g}_0}(\sigma)(z) - B_{\mathbf{g}_0}(\vphi)(z)|\,,
\]
which comes from identifying the tensor $B_{\mathbf{g}_0}(\sigma)(z) - B_{\mathbf{g}_0}(\vphi)(z)$ with the corresponding complex number, as explained in the previous section.
A calculation then shows that the right hand side can be expressed as 
\[
|B_{\mathbf{g}_0}(\sigma)(z) - B_{\mathbf{g}_0}(\vphi)(z)| = |\zeta^2\mathcal{S}f(z)+A(r) -p(r)|\,, \qquad r=|z|, \ \zeta = z/r\,,
\]
where $A(r)$ is defined by \eqref{eq:A(r)}. 

Appealing now to \eqref{eq:|K_g|}, we see that proving \eqref{eq:Bg(sigma-phi)} is equivalent to proving
\[
\left|\zeta^2\mathcal{S}f(z)+A(r)-p(r)\right|+e^{2\sigma(z)}|K(\widetilde{f}(z))| \leq A(r)+p(r) \, .
\]
But in view of the inequality \eqref{eq:p<=A} and the hypothesis \eqref{eq:p-criterion}, we now have
\[
\begin{aligned}
\left|\zeta^2\mathcal{S}f(z)+A(r)-p(r)\right|+e^{2\sigma(z)}|K(\widetilde{f}(z))| 
&\leq |\zeta^2\mathcal{S}f(z)|+|A(r)-p(r)|+e^{2\sigma(z)}|K(\widetilde{f}(z))|\\
&= |\mathcal{S}f(z)|+e^{2\sigma(z)}|K(\widetilde{f}(z)|+A(r)-p(r)\\ 
&\leq A(r)+p(r) \, .
\end{aligned}
\]

\end{proof}

\begin{proof}[Proof of Theorem \ref{prop:convexity}] With 
$u=e^{(\varphi-\sigma)/2}$ and  $v=u^2$ we find using \eqref{eq:Hess-eq} that
\[
{\Hess}_{\mathbf{g}}(v)+vB_{\mathbf{g}}(\sigma-\varphi)=\frac{1}{2}(\Delta_{\mathbf{g}}v)\mathbf{g} \, .
\]
Also $\Delta_{\mathbf{g}}=e^{-2\varphi}\Delta $, where $\Delta$ is the
Euclidean Laplacian, thus
\begin{eqnarray*}
\Delta_{\mathbf{g}}v&=&v\Delta_{\mathbf{g}}(\log v)+\frac{1}{v}||\grad_{\mathbf{g}}v||_{\mathbf{g}}^2\\
&=&ve^{-2\varphi}\Delta(\varphi-\sigma)+\frac{1}{v}||\grad_{\mathbf{g}}v||_{\mathbf{g}}^2\,. 
\end{eqnarray*}
Now using 
$
K_\mathbf{g} = -e^{-2\vphi}\Delta\vphi$, $K = -e^{-2\sigma}\Delta\sigma
$,
and the fact that both curvatures are negative, we can rewrite this as
\[
\Delta_{\mathbf{g}}v=v|K_\mathbf{g}|-ve^{2(\sigma-\varphi)}|K|+\frac{1}{v}||\grad_{\mathbf{g}}v||_{\mathbf{g}}^2 \, ,
\]
and hence
\[
 \Hess_{\mathbf{g}}(v) = -vB_{\mathbf{g}}(\sigma-\varphi)+\frac{v}{2}(|K_\mathbf{g}|-e^{2(\sigma-\varphi)}|K|)\mathbf{g}+
\frac{1}{2v}||\grad_{\mathbf{g}}v||_{\mathbf{g}}^2\,\mathbf{g} \, .
\] 
Therefore, because of the lemma,
\begin{equation}
\Hess_{\mathbf{g}}(v) \geq \left(\frac{v}{2}e^{2(\sigma-\varphi)}|K|
+\frac{1}{2v}||\grad_{\mathbf{g}}v||^2_{\mathbf{g}} \right)\mathbf{g} \, . \label{eq:Hess-v}
\end{equation}
On the other hand, since $v=u^2$
\[
\Hess_{\mathbf{g}}(v)= 2u\Hess_{\mathbf{g}}(u)+2du\otimes du\,,
\] 
and
\[
\frac{1}{2v}||\grad_{\mathbf{g}}v||_{\mathbf{g}}^2 = 2||\grad_{\mathbf{g}}u||_{\mathbf{g}}^2 \, .
\] 
We finally deduce from \eqref{eq:Hess-v} that
\[
2u\,{\Hess}_{\mathbf{g}}(u)+2du\otimes du \geq \frac{v}{2}e^{2(\sigma-\varphi)}|K|\mathbf{g} +
2||\grad_{\mathbf{g}}u||_{\mathbf{g}}^2\,\mathbf{g} \, ,
\]
and since $du \otimes du$ is at most the norm-squared of the gradient,
\[
2u{\Hess}_{\mathbf{g}}(u) \ge \frac{v}{2}e^{2(\sigma-\varphi)}|K|\mathbf{g} \,.
\]
This is equivalent to \eqref{eq:bound-on-hessian}.
\end{proof}

\section{Critical Points, Distortion, and Boundary Behavior}

\label{sec:crit-distortion-bndry}

 We continue to assume that the harmonic mapping $f=h+\bar{g}$ satisfies the univalence criterion \eqref{eq:p-criterion},
\[
|\mathcal{S}f(z)|+e^{2\sigma(z)}|K(\widetilde{f}(z))| \leq 2p(|z|) = \mathcal{S}\Phi(|z|)\,.
\]
We also continue to work with the metric $\mathbf{g} = \Phi'(|z|)^2|dz|^2$ on ${\Bbb D}$. 

In this section we will use the convexity of the function $u$
defined in \eqref{eq:u-for-convexity} in Theorem \ref{prop:convexity} to derive upper bounds on $e^{\sigma(z)} = |h'(z)|+|g'(z)|$. This is the key to obtaining continuous extensions of $f$ and $\widetilde{f}$ to $\partial\Bbb D$ as stated in Theorem \ref{thm:boundary-extension}. The analysis leads to a  more refined understanding of the phenomenon than given by the short, straightforward assertion of the theorem, but distinguishing special cases makes it difficult to collect all the results in a single statement.

\medskip

An important issue is the number of critical points the function $u$
can have, specifically when that number does not exceed one.
 This separates analytic maps, where we can appeal to earlier work, from harmonic maps, where the results of the preceding section will be applied. The distinction is made on the basis of the following lemma.

\begin{lem}If $u$ has more than one critical
point then $\widetilde{f}({\Bbb D})$ is a planar minimal surface.
\end{lem}

\begin{proof} Suppose $z_1$ and $z_2$ are critical points
of $u$. Then, because $u$ is convex, $u(z_1)$ and $u(z_2)$ are
absolute minima, and so is every point on the geodesic segment
$\gamma$ joining $z_1$ and $z_2$ in $\Bbb D$. Hence $
\Hess_{\mathbf{g}}(u)(\gamma', \gamma')=0$, and then Theorem \ref{prop:convexity} implies that
$|K|\equiv 0$ along $\widetilde{\gamma}=\widetilde{f}(\gamma)$. On the
other hand, we know from \eqref{eq:curvature} that
$$ |K| = \frac{4|q'|^2}{|h'|^2(1+|q|^2)^4} \, , $$
hence $q'=0$ along $\gamma$, so $q'$ is identically $0$. This
proves the lemma.
\end{proof}

\medskip
\noindent {\bf Remark.} The proof shows a bit more than stated in the lemma, namely 
that the surface $\Sigma$ will reduce to a plane provided $\Hess_{\mathbf{g}}(u)(\gamma', \gamma')=0$ along any geodesic segment
$\gamma$. In particular, this will be the case if $u$ is constant
along $\gamma$.

\medskip

Thus in the case of multiple critical points the lifted map
$\widetilde{f}$ can be considered as a {\em holomorphic} mapping into a
tilted (complex) plane, and it satisfies the classical Nehari condition
\[
|\mathcal{S}\widetilde{f}(z)| \leq 2p(|z|) \, . 
\]
 A fairly complete treatment of such
mappings, specifically continuous extension
to the boundary, extremal functions (and homeomorphic or
quasiconformal extensions to ${\Bbb C}$) can be found in \cite{gehring:gehring-pommerenke} and
\cite{co:gp}. Briefly, the boundary behavior of $\widetilde{f}$ is of the same character as that of the extremal function $\Phi$ near $x=1$ in the spherical metric (recall that $\Phi(1)=\infty$), a phenomenon we will find to hold  as well when $u$ has one or no critical points. When $p(x) =(1-x^2)^{-2}$ the extremal is a logarithm and $\widetilde{f}$ has a logarithmic modulus of continuity. For all other choices of Nehari functions the extension is H\"older continuous.

\medskip

We next consider the case when $u$ has exactly one critical point. Under that condition, the following lemma is the promised upper bound for $|h'|+|g'|$.

\begin{lem} 
\label{lem:unique-critical-point}
If $u$ has a unique critical point
then there exist constants $a>0$, $b$, and $r_0$\, $(0<r_0<1)$ such that
\begin{equation}
|h'(z)|+|g'(z)| \leq \frac{\Phi'(|z|)}{(a\Phi(|z|)+b)^2} \, ,  \qquad r_0<|z| < 1\, .
\label{eq:bound-for-|h'|+|g'|}
\end{equation}
\end{lem}

\begin{proof} Let $z_0$ be the unique critical point of
$u$. Let $\gamma(s)$ be an arclength parametrized geodesic in the
metric $\mathbf{g}$ starting at $z_0$ in a given direction. Let
$v(s)=u(\gamma(s))$. Because the critical point is unique, it
follows that $v'(s)>0$ for all $s>0$, and hence that there exist an
$s_0>0$ and an $a>0$ such that $v'(s) > a$ for all $s>s_0$.
In turn, $v(s) > as+b$ for some constant $b$ and all
$s>s_0$, and from compactness we can conclude that the constants
$s_0$, $a$, $b$ in this estimate can be chosen uniformly, independent of the
direction of the geodesic starting at $z_0$. In other words,
$$ u(z) \geq ad_{\mathbf{g}}(z,z_0)+b $$
for all $z$ with $d_{\mathbf{g}}(z,z_0)>s_0$, where $d_\mathbf{g}$ denotes distance in the metric $\mathbf{g}$. Then by renaming the
constant $b$, and with a suitable $r_0$, we will have
$$ u(z) \geq ad_{\mathbf{g}}(z,0)+b $$
for all $z$ with $|z| > r_0$. Since $d_{\mathbf{g}}(z,0)=\Phi(|z|)$ the
lemma follows.
\end{proof}

\medskip

We view Lemma \ref{lem:unique-critical-point} as a distortion theorem for harmonic mappings satisfying the univalence criterion, and it is the estimate \eqref{eq:bound-for-|h'|+|g'|} that will allow us to obtain a continuous
extension to the closed disk for the lift $\widetilde{f}$ and
for the harmonic mapping $f$. The modulus of continuity of the
extension depends on particular properties of the function $\Phi$, and 
ranges from a logarithmic modulus of continuity, to H\"older and to
Lipschitz continuity. 

We begin by observing that since the
function $(1-x^2)^2p(x)$ is positive and nonincreasing on $[0,1)$ we can form the limit
\[
\lambda= \lim_{x\rightarrow 1^-}(1-x^2)^2p(x)\,.
\]
It was shown
in \cite{co:noncomplete} that $\lambda\leq 1$ and that $\lambda=1$ if and only if
$p(x)=(1-x^2)^{-2}$. Consider first this case, when 
$p(x)=(1-x^2)^{-2}$. Then
$$
\Phi(x)=\frac{1}{2}\log\frac{1+x}{1-x}\,,
$$
and \eqref{eq:bound-for-|h'|+|g'|} amounts to
\begin{equation}
|h'(z)|+|g'(z)| \leq
\frac{1}{(1-|z|^2)\left(\displaystyle{\frac{a}{2}\log\frac{1+|z|}{1-|z|}}+b\right)^2}
\,, \qquad   r_0 <|z| < 1\, . 
\label{eq:gp-argument}
\end{equation}

 From here, to show that $\widetilde{f}$ extends continuously to the closed disk, we follow the argument in \cite{gehring:gehring-pommerenke} and integrate along a hyperbolic geodesic. Let
 $\varrho$ be the hyperbolic segment
joining two points $z_1$ and $z_2$ in $\Bbb D$. Then $\varrho$ has Euclidean length
$\ell\leq\frac{\pi}{2}|z_1-z_2|$ and
$\min(s,\ell-s)\leq\frac{\pi}{2}(1-|z|)$ for each $z\in\varrho$,
where $s$ is the Euclidean arclength of the part of $\varrho$
between $z_1$ and $z$. Suppose that $z_1$ and $z_2$ are such that $\varrho$ is contained in the
annulus $r_0<|z|<1$. The distance $|\widetilde{f}(z_1)-\widetilde{f}(z_2)|$ in ${\Bbb R}^3$ is less than the metric distance between $\widetilde{f}(z_1)$ and $\widetilde{f}(z_2)$ on the surface $\Sigma$, and this in turn is less than the length of $\varrho$ in the metric $(|h'|+|g'|)|dz|$ on $\Bbb D$. Thus we may use \eqref{eq:bound-for-|h'|+|g'|} and \eqref{eq:gp-argument} and write
\begin{eqnarray*}
|\widetilde{f}(z_1)-\widetilde{f}(z_2)| &\leq&
\int_{\varrho}(|h'(z)|+|g'(z)|)|dz| \leq \int_{\varrho}
\frac{|dz|}{(1-|z|^2)({\frac{a}{2}\log\frac{1+|z|}{1-|z|}}+b)^2}\\
&\leq& C\int_0^{\ell/2}\frac{ds}{(1-s)({\frac{a}{2}\log\frac{1}{1-s}}+b)^2}\,,
\end{eqnarray*}
for some constant $C$ independent of $z_1$ and $z_2$. 

Integration, together with the bound $\ell\leq\frac{\pi}{2}|z_1-z_2|$,
yields
\[
|\widetilde{f}(z_1)-\widetilde{f}(z_2)| \leq
\frac{C'}{\log\displaystyle{\frac{1}{|z_1-z_2|}}} \, , 
\]
for another constant $C'$.  This implies that $\widetilde{f}$ has an extension to $\overline{\Bbb D}$ that is uniformly
continuous. Since 
\begin{equation}
|f(z_1)-f(z_2)| \leq
|\widetilde{f}(z_1)-\widetilde{f}(z_2)|\,, \label{eq:f-continuous}
\end{equation}
the same is true for the harmonic mapping $f$. In all the cases that we consider, it is simply the inequality \eqref{eq:f-continuous} that is used to obtain a continuous extension for $f$ from one for $\widetilde{f}$. 

\medskip

Suppose now that $\lambda<1$. It was shown in \cite{co:noncomplete} that
\begin{equation}
\lim_{x\rightarrow 1^-}(1-x^2)\frac{\Phi''(x)}{\Phi'(x)} = 2(1+\sqrt{1-\lambda})=2\mu \, . 
\label{eq:2mu}
\end{equation}
Note that $1<\mu\leq 2$. It follows from \eqref{eq:2mu} that for any
$\varepsilon>0$ there exists $x_0\in(0,1)$ such that
\[
\frac{\mu-\varepsilon}{1-x} \leq \frac{\Phi''(x)}{\Phi'(x)}\leq \frac{\mu+\varepsilon}{1-x}\,, \qquad x> x_0 \, .
\]
This implies that
\[
\frac{1}{(1-x)^{\mu-\varepsilon}} \leq \Phi'(x) \leq \frac{1}{(1-x)^{\mu+\varepsilon}}\, , \qquad x> x_0 \, ,
\]
so that
\begin{equation}
\frac{\Phi'(x)}{(a\Phi(x)+b)^2} \leq \frac{C}{(1-x)^{\alpha+3\varepsilon}} \, , 
\label{eq:F'-upper}
\end{equation}
where $\alpha=2-\mu=1-\sqrt{1-\lambda}$ and $C$ depends on $a$, $b$
and the values of $\Phi$ at $x_0$. The estimate in \eqref{eq:F'-upper} together
with the technique of integrating along a hyperbolic
segment implies now that
\[
|\widetilde{f}(z_1)-\widetilde{f}(z_2)| \leq C|z_1-z_2|^{1-\alpha-3\varepsilon}
= C|z_1-z_2|^{\sqrt{1-\lambda}-3\varepsilon} \, ,
\] for all points
$z_1,z_2$ whose joining hyperbolic geodesic segment
is contained in the annulus $\max\{r_0,x_0\}<|z|<1$. This shows
that $\widetilde{f}$, and hence $f$, admits a continuous extension to the
closed disk, with at least a H\"older modulus of continuity.

Observe that the left hand side of \eqref{eq:F'-upper} is the derivative of a M\"obius transformation of $\Phi$. Thus the modulus of continuity of $\widetilde{f}$, as derived from this bound, is, in essence, that of the extremal at $x=1$ in the spherical metric.

We also remark that if additionally it is known that $x=1$ is a
regular singular point of the differential equation $u''+pu=0$,
then from the analysis of the Frobenius solutions at
$x=1$ it follows that
\[
\Phi'(x) \sim \frac{1}{(1-x)^{\mu}}\,, \qquad  x\rightarrow 1 \, . 
\]
This then provides for a H\"older continuous extension with the `best' exponent when
$\lambda >0$, and a Lipschitz continuous extension when $\lambda=0$.

\medskip

Recall that the preceding arguments were carried out under the umbrella of Lemma \ref{lem:unique-critical-point}, when $u$ has a unique critical point. In this case the extension to the closed disk is continuous (or more) in the Euclidean metric. We finally treat the case when $u$ has no critical points, and here the situation is somewhat different. 

A fairly straightforward argument gives half a result, so to speak. Suppose that $u$ has no critical points and let
$u_{\theta}(s)=u(r(s)e^{i\theta})$, where $r(s)$ is the arclength
parametrization of the radius $[0,1)$ in the metric $\mathbf{g}$. By assumption, the gradient
of $u$ at the origin does not vanish, and hence $u_{\theta}'(0)>0$ for
all arguments $\theta$ in an open half-plane; to be specific, say $u_{\theta}'(0)>0$ for $0<\theta<\pi$. Furthermore, because $u$ cannot be
constant along any geodesic unless the surface $\Sigma$ reduces to
a plane, it follows that $u_{\theta}'(s)>0$ for $s>0$ and
$\theta=0$ or $\theta =\pi$. Now again by compactness we see that
there is an $s_0$ such that $u_{\theta}'(s)>a>0$ for $s\ge s_0$ and all $0\leq \theta \leq
\pi$. From here we can pick up the proof of Lemma \ref{lem:unique-critical-point} and deduce that \eqref{eq:bound-for-|h'|+|g'|}
holds for all $|z|>r_0$ with ${\rm Im}\{z\} \geq 0$. To reiterate, this then provides a continuous extension of $\widetilde{f}$ and of $f$ to the upper half of $\overline{{\Bbb D}}$.

To get beyond this half-disk result we will show that for any
radius $[0,e^{i\theta_0})$ there
exists a M\"obius transformation $T$ of ${\Bbb R}^3$ such that the
conformal factor  associated with the conformal lift
$T\circ\widetilde{f}$ satisfies a version of \eqref{eq:bound-for-|h'|+|g'|} in an angular sector
$|\theta-\theta_0|<\delta$ about the radius -- why we need the extra M\"obius transformation will emerge presently.   Thus $T\circ\widetilde{f}$ will exhibit
the appropriate continuous extension. It is because we have to allow for a shift of the range by a M\"obius transformation of ${\Bbb R}^3$ that the summary result on boundary extensions, Theorem \ref{thm:boundary-extension}, is stated  to assert that $\widetilde{f}$ and $f$ have extensions to the closed disk that are continuous in the {\em spherical} metric.

\medskip

Let $T$ be a M\"obius transformation of ${\Bbb R}^3$. Since $T\circ \widetilde{f}$, though conformal, may not be the lift of a harmonic mapping we do not have the basic convexity result Theorem \ref{prop:convexity} as a starting point. We shall first find a substitute that holds along a radius.  

Let $e^{\tau}$ be the conformal factor associated with $T\circ \widetilde{f}$, that is, $e^{\tau(z)}=|DT(\widetilde{f}(z))|e^{\sigma(z)}$, and for fixed $\theta$ let $\tau_{\theta}(r)=\tau(re^{i\theta})$. Along each radius we have:

\begin{lem}
\label{lem:substitute-convexity}
Let $r=r(s)$ be the arclength
parametrization of $[0,1)$ in the metric $\mathbf{g}=\Phi'(|z|)^2|dz|^2$. Then the function
\begin{equation}
\omega(s)=\sqrt{\frac{\Phi'(r(s))}{e^{\tau_{\theta}(r(s))}}}  \label{eq:omega-is-convex}
\end{equation}
is a convex function of $s$.
\end{lem}

Before giving the proof we note that it is easy to identify the arclength parametrization of a radius in the metric $\mathbf{g}$. Since $\Phi$ is radial the length $s$ of the radius $[0,re^{i\theta}]$ is 
\[
s=\int_0^r \Phi'(t)\,dt = \Phi(r)\,.
\]
Thus if $\Psi$ denotes the inverse of $\Phi$, then $r(s) = \Psi(s)$ in the notation of the lemma.

\begin{proof} Let  $w(r)=\Phi'(r)^{-1/2}$,  so that $w''+pw=0$, and let $g(r)=(T\circ \widetilde{f})(re^{i\theta})$ and note that 
$|g'(r)|=e^{\tau_{\theta}(r)}$. A straightforward calculation
shows that the function $v(r)=|g'(r)|^{-1/2}$ satisfies the
equation $ v''+qv=0 $ with
\[
q = \frac{1}{2}S_1g - \frac{1}{4}\left(\frac{|g''|^2}{|g'|^2}-
\frac{\langle g'', g'\rangle^2}{|g'|^4}\right)\,,
\]
where $S_1$ is Ahlfors' Schwarzian \eqref{eq:S1}. The quantity in parentheses is nonnegative, and  by M\"obius invariance, $S_1g(r) = S_1\widetilde{f}(re^{i\theta})$. Since $S_1\widetilde{f}(re^{i\theta})\le 2p(re^{i\theta})$ by assumption, we conclude that $q\le p$. 

Introducing $\Psi=\Phi^{-1}$ as above, we write
\[
\omega(s)=\sqrt{\frac{\Phi'(r(s))}{e^{\tau_{\theta}(r(s))}}}=\frac{v(\Psi(s))}{w(\Psi(s))}\,,
\]
and now one can check that
\[
\omega'' = (p-q)w^4 \omega\,,
\]
where differentiation is with respect to $s$ and the quantity $(p-q)w^4$ is to be evaluated at $\Psi(s)$.  Since
$p\geq q$, it follows that $\omega$ is a convex function of $s$.
\end{proof}

To obtain an estimate of the type \eqref{eq:bound-for-|h'|+|g'|} along a given radius
$[0,e^{i\theta})$ it will suffice to show that the function $\omega$ in
\eqref{eq:omega-is-convex} has $\omega'(0) > 0$. Indeed, if this derivative is positive
then, once more by convexity, $\omega(s) \geq as+b$ for some constant $b$ and some positive
constant $a$. With this,
\begin{equation}
e^{\tau_{\theta}(r)} \leq \frac{\Phi'(r)}{(a\Phi(r)+b)^2} \,, \label{eq:bound-for-tau}
\end{equation}
corresponding to \eqref{eq:bound-for-|h'|+|g'|}. Then, if for a given angle $\theta_0$ the derivative of $\omega$ at
zero is positive, by continuity it will remain positive for
all angles $\theta$ close to $\theta_0$. Therefore \eqref{eq:bound-for-tau} will hold
uniformly in an angular sector, leading to a continuous extension of $\widetilde{f}$ and $f$ to the part of the boundary within the sector, just as before. 

To finish the argument, we thus need to
find a  M\"obius transformation $T$ of ${\Bbb R}^3$ so that,
for a given angle $\theta_0$, the function $\omega$ has positive
derivative at zero. Using $\Phi''(0)=0$ we are therefore required to have
$\tau_{\theta_0}'(0)<0$. After a translation, a rotation and a
dilation, we can assume that the curve
$\widetilde{f}_{\theta_0}(r)=\widetilde{f}(re^{i\theta_0})$ has
$\widetilde{f}_{\theta_0}(0)=0$, $\widetilde{f}_{\theta_0}'(0)=(1,0,0)$
and $\widetilde{f}_{\theta_0}''(0)=(\alpha, \beta, 0)$. Let $T$ be the
extension as a M\"obius transformation of ${\Bbb R}^3$ of the complex M\"obius
map 
\[
z\mapsto \displaystyle{\frac{z}{1+cz}}\,,
\] 
where $c=(1+\alpha)/2$ and we identify $z=x+iy$ with the point
$(x,y,0)$. Since up to order 2 the curve $\widetilde{f}_{\theta_0}(r)$
lies in the $(x,y)$-plane, for the purposes of our calculations, which involve derivatives of order 2 at most, we may replace the
curve $(T\circ\widetilde{f}_{\theta_0})(r)$ with the curve
\[
w(r)=\displaystyle{\frac{z(r)}{1+cz(r)}}\,,
\] 
where $z(r)$ satisfies
$z(0)=0, z'(0)=1, z''(0)=\alpha+i\beta$. Up to an error of
order $O(r^2)$ we have that $\tau_{\theta_0}(r) = \log
|w'(r)|=\log |z'(r)| -2\log|1+cz(r)|$, from which
$\tau_{\theta_0}'(0)= \alpha-2c=-1<0$, as desired.

\smallskip
\noindent {\bf Remark.} The M\"obius transformation $T$ may send
some point on the surface $\Sigma$ to the point at infinity, but
such a point cannot lie in the image of $\widetilde{f}_{\theta_0}(r)$.
Indeed, once we have ensured that the function $\omega$ in Lemma \ref{lem:substitute-convexity}
has $\omega'(0)>0$, then the estimate \eqref{eq:bound-for-tau} will imply that the
curve $(T\circ\widetilde{f}_{\theta_0})(r)$ has finite length, and it is therefore
impossible for it to reach the point at infinity.

\medskip

The arguments in this section, supported by the results of the previous two sections, have proved Theorem 2 as stated, that $f$ and $\widetilde{f}$ have spherically continuous extensions to $\overline{\Bbb D}$. Much more detail has been obtained {\em en route}, and we conclude with an expanded, if admittedly underspecified version of the theorem that we hope helpfully captures the main points.

\begin{thm2'}
Suppose $f$ satisfies the univalence criterion 
\[
|\mathcal{S}f(z)|+e^{2\sigma(z)}|K(\widetilde{f}(z))| \leq 2p(|z|) = \mathcal{S}\Phi(|z|)\,,
\]
 with extremal function $\Phi$, and let $\lambda = \lim_{x\rightarrow 1^-}(1-x^2)^2p(x)$. Then $f$ and $\widetilde{f}$ have extensions to $\overline{{\Bbb D}}$ that are continuous with respect to the spherical metric. 
The modulus of continuity of each is of the same type as that of $\Phi(x)$ near $x=1$ in the spherical metric. If $\lambda=1$ it is logarithmic. If $\lambda <1$ it is H\"older with an exponent that depends on $\lambda$.
\end{thm2'}

\section{The Catenoid and Extremal Lifts}
\label{section:examples}

The principal work of this section is to consider some examples that show our
results are sharp. The constructions are based on what we know from the analytic case, and on one of the earliest minimal surfaces to be studied, the catenoid.  According to Theorem 3, the minimal surface corresponding to the lift of an extremal mapping must contain a Euclidean circle or line as a circle of curvature. Catenoids enter naturally into the discussion of extremal lifts because they are the \emph{unique} minimal surfaces containing a Euclidean circle as a line of curvature, as we shall now show.

\begin{lem} \label{lemma:catenoid}
If a minimal surface $\Sigma$ contains a part of a Euclidean circle or line as a line of curvature, then $\Sigma$ is contained in a catenoid or a plane.
\end{lem}

\begin{proof}
>From the theory of minimal surfaces we will need a uniqueness result associated with the \emph{Bj\"orling problem} of finding a minimal surface with a prescribed normal strip. This may be stated as follows. Let $I\subset\mathbb{R}$ be an open interval. A
real-analytic strip $S=\{({\bold c}(t),{\bold n}(t)):t\in I\}$ in $\mathbb{R}^3$ consists of a
real-analytic curve ${\bold c}:I\rightarrow\mathbb{R}^3$ with ${\bold c}'(t)\neq
0$ and a real-analytic vector field ${\bold n}:I\rightarrow\mathbb{R}^3$
along ${\bold c}$, with $|{\bold n}(t)|\equiv 1$ and $\langle {\bold c}'(t),{\bold n}(t)\rangle
\equiv 0$. The problem is to find a
parametrized minimal surface ${\bold X}:\Omega\rightarrow\mathbb{R}^3$
with $I\subset\Omega \subset \mathbb{R}^2$, such that ${\bold X}(x,0)={\bold c}(x),\, {\bold N}(x,0)={\bold n}(x)$ for
$x\in I$, where ${\bold N}(x,y)$ is a unit normal vector field along the surface. The result we need is that for real-analytic data,
Bj\"{o}rling's problem admits exactly one solution; see \cite{dhkw:minimal}, p. 121, where the solution is expressed in closed form in terms of the data defining the strip.

Suppose now that $C$ is a Euclidean circle, part (or all) of which is a line of curvature of a
minimal surface $\Sigma$. Let $\Pi$ be the plane containing $C$, and
let ${\bold N}$, ${\bold n}_0$ be, respectively, unit normal vectors to $\Sigma$
and to $\Pi$. It follows from the classical Theorem of Joachimstahl (see, for example, \cite{doCarmo:diff-geom} p. 152) that the normal vector of a surface along a planar
curve which is a line of curvature forms a constant angle with the
normal to the plane of the curve. Therefore, in our case, ${\bold N}$ and ${\bold n}_0$ form a
constant angle along $C$, say $\alpha$.

Let $\Sigma_0$ be a catenoid and consider the situation above for $\Sigma_0$. All of the circles of revolution on $\Sigma_0$ are lines of curvature for $\Sigma_0$, and the angle between the normal to $\Sigma_0$ and the plane of any such circle decreases from $\pi/2$, for the circle around the waist of $\Sigma_0$, down toward $0$ as the circles move out toward infinity. Choose one circle, $C_0$, where the angle is $\alpha$. Via Euclidean similarities we can assume that $C_0$ and $C$ coincide, and so both surfaces $\Sigma$ and $\Sigma_0$ have the same unit normal vector fields along $C=C_0$. 

Next let ${\bold X}$ and ${\bold X}_0$ be conformal parametrization of $\Sigma$  and $\Sigma_0$, respectively, covering a part of the circle $C =C_0$, say $C'$. We arrange the parametrization of $\Sigma_0$ so that ${\bold X}_0^{-1}(C') =I$, an open interval. Because the preimage $c'={\bold X}^{-1}(C')$ is a real analytic, simple curve, there is an invertible holomorphic map, $h$, of a neighborhood of $I$ to a neighborhood of $c'$ with ${\bold X}(h(t))={\bold X}_0(t)$, $t \in I$. It now follows from the uniqueness of the solution to Bj\"orling's problem that ${\bold X} = {\bold X}_0\circ h^{-1}$ on a neighborhood of $c'$. Thus $\Sigma$ and $\Sigma_0$ coincide near $C'$, and hence $\Sigma$ is a portion of the catenoid. If $\Sigma$ contains part of a Euclidean line as a line of curvature, instead of a circle, then a similar argument shows that $\Sigma$ must be contained in a plane.  

\end{proof}

With this result, we now have the following corollary of Theorem \ref{thm:extremal} on extremal maps.

\begin{cor}
Let $f$ be an extremal mapping. Then $\tilde{f}(\mathbb{D})$ is contained in a catenoid or a plane.
\end{cor}

We now proceed with examples. The case of an extremal lift mapping the disk into a plane is essentially the case of an analytic extremal, and examples there have been studied. More interesting for the present paper are extremal mappings into a catenoid, where the analytic case can still serve as a guide.

\medskip

\noindent {\bf Example 1:} The choice $p(x)=\pi^2/4$ gives Nehari's univalence criterion \eqref{eq:nehari-pi^2/2} for analytic functions in the disk. In the analytic case an extremal mapping is 
\[
F(z)=\displaystyle{\frac{2}{\pi}\tan(\frac{\pi}{2}z)}\,,
\]
which maps ${\Bbb D}$ to a horizontal strip. For harmonic maps the criterion \eqref{eq:p-criterion} becomes
\begin{equation}
|\mathcal{S}f(z)|+e^{2\sigma(z)}|K(\widetilde{f}(z))| \leq \frac{\pi^2}{2} \,. 
\label{eq:pi^2/2-criterion}
\end{equation}
To show that the criterion is sharp we will work with the harmonic mapping
\[
f(z)=h(z)+\overline{g(z)}=ce^{\pi z}+\frac{1}{c}e^{-\pi\bar{z}} \, , 
\]
for a positive constant $c$ to be chosen later. The lift $\widetilde{f}$ maps ${\Bbb D}$ into the catenoid parametrized by
\begin{eqnarray*}
U(x,y) &=& (ce^{\pi x} + \frac{1}{c}e^{-\pi x})\cos \pi y \\
V(x,y) &=& (ce^{\pi x} + \frac{1}{c}e^{-\pi x})\sin \pi y\\
W(x,y) &=& 2\pi x
\end{eqnarray*}
with $z=x+iy$. The lift fails to be univalent at $\pm i$, with $\widetilde{f}(i) = \widetilde{f}(-i) = (-(c+\frac{1}{c}),0,0)$. In fact, the diameter $-1 \le y \le 1$ maps to the circle $U^2+V^2 =(c+\frac{1}{c})^2$, $W=0$ on the surface. This is one of the circles of revolution of the catenoid, and it is a line of curvature as guaranteed by Theorem 3.

To see what happens with \eqref{eq:pi^2/2-criterion}, we find first that
\[
e^{\sigma(z)}=|h'(z)|+|g'(z)|=\pi(ce^{\pi x}+\frac{1}{c}e^{-\pi x}) \,,
\]
and then for the Schwarzian, 
\[
\mathcal{S}f(z)= -\frac{\pi^2}{2}+4\pi^4e^{-2\sigma(z)} \,.
\]
For the curvature term,
\[
e^{2\sigma(z)}|K(\widetilde{f}(z)|=\Delta\sigma(z)=4\pi^4e^{-2\sigma(z)} \,. 
\]
Therefore, \eqref{eq:pi^2/2-criterion} will be satisfied provided
\[
\left|-\frac{\pi^2}{2}+4\pi^4e^{-2\sigma(z)}\right|+4\pi^4e^{-2\sigma(z)} \leq \frac{\pi^2}{2} \, . 
\]
This will be the case if $c>(1+\sqrt{2})e^\pi=55.866\dots$, because then
\[
\left|-\frac{\pi^2}{2}+4\pi^4e^{-2\sigma(z)}\right| = \frac{\pi^2}{2}-4\pi^4e^{-2\sigma(z)} \, , 
\]
and for $c$ in this range
\[
|\mathcal{S}f(z)|+e^{2\sigma(z)}|K(\widetilde{f}(z))| \equiv \frac{\pi^2}{2}\,.
\]

By modifying this construction slightly we can also show that the constant $\pi^2/2$ is best possible. For the harmonic mapping take
\[
f(z)=ce^{t\pi z}+\frac1c e^{-t\pi \bar{z}}\,,\quad t >0\,.
\]
 Then 
\[
e^{\sigma(z)}=t\pi(ce^{t\pi x}+\frac{1}{c} e^{-t\pi x}) \, ,
\]
and
\[
\mathcal{S}f=-t^2\frac{\pi^2}{2}+4t^4\pi^4e^{-2\sigma}\,,
\]
while
\[
e^{-2\sigma}|K|=\Delta\sigma=4t^4\pi^4e^{-2\sigma}\, . 
\]
Therefore, if $c>56$,
\[
|\mathcal{S}f|+e^{-2\sigma}|K|=\left|-t^2\frac{\pi^2}{2}+4t^4\pi^4e^{-2\sigma}\right|+
4t^4\pi^4e^{-2\sigma} = t^2\frac{\pi^2}{2}\,.
\]
But as soon as $t>1$ both the maps $f$ and $\widetilde{f}$ fail to be
univalent in ${\Bbb D}$.

\medskip
\noindent {\bf Example 2:} Portions of the 
catenoid also provide examples for other Nehari functions.  
We discuss a general procedure. Let $p$ be a Nehari function 
that is the restriction to $(-1,1)$ of
an analytic function $p(z)$ in the disk  with the property  $|p(z)|\leq
p(|z|)$. Typical examples are $p(z)=(1-z^2)^{-2}$ and
$p(z)=2(1-z^2)^{-1}$. The extremal map in such a case, say $F$, can be normalized in the same way as the extremal $\Phi$ in \eqref{eq:Phi}, and is analytic, odd, and
univalent in the disk and satisfies $SF(z)=2p(z)$. The image
$F({\Bbb D})$ is a ``parallel strip'' like domain, symmetric with respect
to both axis, and containing the entire real line; see \cite{co:gp}. Let
\[
G(z)=\frac{cF(z)+i}{cF(z)-i} \, , 
\]
where $c>0$ is to be chosen later and sufficiently small so that
$i/c \notin F({\Bbb D})$ (it can be shown that the map $F$ is always
bounded along the imaginary axis). The function $G$ maps ${\Bbb D}$ onto
a simply-connected domain containing the unit circle minus the
point 1. Let
\[
f(z)=h(z)+\overline{g(z)}= G(z)+\overline{G(-z)} \, . 
\label{eq:extremal-G}
\]
Since $F$ is odd, 
\[
f(z)=G(z)+\frac{1}{\overline{G(z)}}\,,
\]
and it follows that
the lift $\widetilde{f}$ parametrizes the catenoid with the unit
circle $|G|=1$ mapped onto the circle of symmetry of the catenoid.
We also have
\[
\begin{aligned}
e^{\sigma}=|G'(z)|+|G'(-z)|&=\frac{2c|F'(z)|}{|cF(z)-i|^2}+\frac{|2cF'(-z)|}{|cF(-z)-i|^2}\\
&=\frac{2c|F'(z)|}{|cF(z)-i|^2}+\frac{2c|F'(z)|}{|cF(z)+i|^2} \, ,
\end{aligned}
\]
using again that $F$ is odd. A somewhat tedious calculation shows that
\[
\mathcal{S}f=2(\sigma_{zz}-\sigma_z)=\mathcal{S}F-\frac{4c^2(1+c^2\overline{F}^2)(F')^2}{(1+c^2F^2)(1+c^2|F|^2)^2}
\, , 
\]
and
\[
e^{2\sigma}|K|=\frac{4c^2|F'|^2}{(1+c^2|F|^2)^2} \, . 
\]
Condition \eqref{eq:p-criterion} now reads
\begin{equation}
\left|\mathcal{S}F(z)-\frac{4c^2(1+c^2\overline{F(z)}^2)F'(z)^2}{(1+c^2F(z)^2)(1+c^2|F(z)|^2)^2}\right|
+\frac{4c^2|F'(z)|^2}{(1+c^2|F(z)|^2)^2} \leq \mathcal{S}F(|z|) \, . 
\label{eq:p-criterion-F}
\end{equation}

Suppose, for example, we let $p(z)=(1-z^2)^{-2}$, for
which the extremal function is 
\[
F(z)=\frac{1}{2}\log\frac{1+z}{1-z}\,.
\]
One has 
\[
F'(z) = \frac{1}{1-z^2}\,,\quad \mathcal{S}F(z) = \frac{2}{(1-z^2)^2}\,,
\]
and \eqref{eq:p-criterion-F} becomes
\[
\begin{aligned}
\left|\frac{2}{(1-z^2)^2}-\frac{4c^2(1+c^2\overline{F(z)}^2)}{(1-z^2)^2(1+c^2F(z)^2)(1+c^2|F(z)|^2)^2}\right|
&+
\frac{4c^2}{|1-z^2|^2(1+c^2|F(z)|^2)^2}\\
 \leq  \frac{2}{(1-|z|^2)^2}\, ,
\end{aligned}
\]
 which reduces to
\begin{equation}
\left|1-\frac{2c^2(1+c^2\overline{F(z)}^2)}{(1+c^2F(z)^2)(1+c^2|F(z)|^2)^2}\right|
+\frac{2c^2}{(1+c^2|F(z)|^2)^2} \leq \frac{|1-z^2|^2}{(1-|z|^2)^2} \,
. 
\label{eq:p-criterion-F-reduced}
\end{equation} 
We comment at once that equality holds here if $z$ is real and if $c$ is sufficiently small, for  both sides of the inequality are then just $1$. The task is to show that  \eqref{eq:p-criterion-F-reduced} holds for all $z \in {\Bbb D}$.

Let
\[
\zeta=\displaystyle{\frac{2c^2(1+c^2\overline{F(z)}^2)}{(1+c^2F(z)^2)(1+c^2|F(z)|^2)^2}}\,.
\]
We establish the following estimates.

\begin{lem} If $c$ is small then there exist
absolute constants $A,B,C$ such that
\begin{equation}
|1-{\rm Re}\{\zeta\}| \leq 1-\frac{2c^2}{(1+c^2|F(z)|^2)^2} +Ac^4|{\rm Im}\{F(z)\}|^2 \, , 
\label{eq:estimate-for-Re}
\end{equation}
\begin{equation}
|{\rm Im}\{\zeta\}| \leq Bc^3|{\rm Im}\{F(z)\}| \, , 
\label{eq:estimate-for-Im}
\end{equation}
and
\begin{equation}
|1-\zeta| \leq  1-\frac{2c^2}{(1+c^2|F(z)|^2)^2} +Cc^4|{\rm Im}\{F(z)\}|^2 \, . 
\label{eq:estimate-for-|1-zeta|}
\end{equation}
\end{lem}

\begin{proof} We write
\begin{eqnarray*}
&&1-{\rm Re}\{\zeta\}=1-\frac{2c^2}{(1+c^2|F(z)|^2)^2}{\rm Re}\left\{\frac{1+c^2\overline{F(z)}^2}{1+c^2F(z)^2}\right\}
\\ 
&& \hskip,25in =
1-\frac{2c^2}{(1+c^2|F(z)|^2)^2}+
\frac{2c^2}{(1+c^2|F(z)|^2)^2}\left(1-{\rm
Re}\left\{\frac{1+c^2\overline{F(z)}^2}{1+c^2F(z)^2}\right\}\right) \,,
\end{eqnarray*} 
which after some calculations gives
\[
1-{\rm Re}\{\zeta\} = 1-\frac{2c^2}{(1+c^2|F(z)|^2)^2}+
\frac{2c^4|cF(z)+c\overline{F(z)}|^2}{|1+c^2F(z)^2|^2(1+c^2|F(z)|^2)^2}|F(z)-\overline{F(z)}|^2
\, . 
\]
The inequality \eqref{eq:estimate-for-Re} follows from this because the quantity
\[
\displaystyle{\frac{|cF(z)+c\overline{F(z)}|^2}{|1+c^2F(z)^2|^2(1+c^2|F(z)|^2)^2}}
\]
is uniformly bounded for small $c$.

To establish \eqref{eq:estimate-for-Im} we have
\begin{eqnarray*}
{\rm Im}\{\zeta\} &=& \frac{2c^2}{(1+|F(z)|^2)^2}{\rm Im}\left\{\frac{1+c^2\overline{F(z)}^2}{1+c^2F(z)^2}\right\}\\
&=&
2c^3\frac{(2+c^2F(z)^2+c^2\overline{F(z)}^2)(cF(z)+c\overline{F(z)})}{|1+c^2F(z)^2|^2(1+|F(z)|^2)^2}
\frac{\overline{F(z)}-F(z)}{2i} \, , 
\end{eqnarray*}
from which \eqref{eq:estimate-for-Im} follows since
\[
\displaystyle{\frac{(2+c^2F(z)^2+c^2\overline{F(z)}^2)(cF(z)+c\overline{F(z)})}{|1+c^2F(z)^2|^2(1+|F(z)|^2)^2}}
\]
is uniformly bounded for small $c$. Finally, \eqref{eq:estimate-for-|1-zeta|} is a
consequence of \eqref{eq:estimate-for-Re} and \eqref{eq:estimate-for-Im} because for $|\zeta|$ small,  $\zeta=x+iy$,  we have
$|1-\zeta| \leq |1-x|+2y^2$.
\end{proof}

With this lemma we can now obtain \eqref{eq:p-criterion-F-reduced} via an analysis
along the level sets of the function $|1-z^2|/(1-|z|^2)$.
The set of points $z$ where
$|1-z^2|/(1-|z|^2)=\sqrt{1+t^2}$, $t>0$, corresponds to a
pair of arcs of circles through $\pm 1$ centered at $\pm i/t$ with
radius $\sqrt{1+t^2}/t$. Because of the symmetry of \eqref{eq:p-criterion-F-reduced} it
suffices to consider the part of the upper arc, call it $\gamma$.
The arc $\gamma$ intersects the imaginary axis at
$it/(1+\sqrt{1+t^2})$ and is mapped under $F$ to the horizontal
line ${\rm Im}\{F\}=\tan^{-1}(s)$, where $s=t/(1+\sqrt{1+t^2})
\leq t$. Thus ${\rm Im}\{F\} \leq \tan^{-1}t \leq t$. From
\eqref{eq:estimate-for-|1-zeta|} it follows that along $\gamma$
$$ |1-\zeta|  \leq  1-\frac{2c^2}{(1+c^2|F(z)|^2)^2} +Cc^3t^2 \, ,$$ so that the
left hand side of \eqref{eq:p-criterion-F-reduced} is bounded above by $1+Cc^4t^2<1+t^2$ for $c$  sufficiently small.

We have shown that the criterion
\[
|\mathcal{S}f(z)|+e^{2\sigma(z)}|K(\widetilde{f}(z))| \leq \frac{2}{(1-|z|^2)^2} 
\]
is sharp. By adapting an example given by Hille \cite{hille:remark} (which accompanied Nehari's original paper) we can also show that the constant $2$ in the numerator of the right hand side is best possible. For this take
\[
F(z)=\left(\frac{1+z}{1-z}\right)^{i\varepsilon} \, ,
\]
which is  far from univalent in the unit disk if $\varepsilon>0$; the value $1$ is assumed infinitely often. 
Note that
\[
F'(z)=\frac{2i\varepsilon}{1-z^2}F(z)\,, \quad  \mathcal{S}F(z) = \frac{2(1+\varepsilon^2)}{(1-z^2)^2} \, , \]
and 
\[
e^{-\varepsilon\frac{\pi}{2}}\leq |F(z)|\leq
e^{\varepsilon\frac{\pi}{2}} \, . 
\]
>From these it is easy to see that equation \eqref{eq:p-criterion-F} will
be satisfied if the right hand side is replaced by
\[
\frac{2+\delta}{(1-|z|^2)^{2}}\,,
\]
 where we can make $\delta>0$
arbitrarily small if $c$ and $\varepsilon$ are each sufficiently
small.
 
\smallskip

For one final example, if we take $p(z)=(1-z^2)^{-1}$ then an extremal map is
\[
F(z)=\displaystyle{\int_0^z\frac{d\zeta}{(1-\zeta^2)^2}} \,,
\]
and similar calculations show that \eqref{eq:p-criterion-F} will be satisfied for sufficiently small $c$ with equality for $z$ real.

\bibliography{lift}

\end{document}